\documentclass[11pt]{article}
\title
{The phase transition in random graphs -- a simple proof}
\author{Michael Krivelevich
\thanks{School of Mathematical Sciences, Raymond and Beverly
Sackler Faculty of Exact Sciences, Tel Aviv University, Tel Aviv,
69978, Israel. Email: krivelev@post.tau.ac.il. Research supported in
part by a USA-Israel BSF grant and by a grant from the Israel
Science Foundation.}
\and Benny Sudakov
\thanks{Department of Mathematics, UCLA, Los Angeles, CA 90095.
Email: bsudakov@math.ucla.edu. Research supported in part by NSF
grant DMS-1101185, by AFOSR MURI grant FA9550-10-1-0569  and by a
USA-Israel BSF grant.}
 }

\usepackage{amsmath, amsfonts}

\begin{document}
\bibliographystyle{plain}
\maketitle
\newtheorem{thm}{Theorem}
\newtheorem{defin}{Definition}
\newtheorem{lemma}{Lemma}
\newtheorem{thmtool}{Theorem}[section]
\newtheorem{corollary}[thmtool]{Corollary}
\newtheorem{lem}[thmtool]{Lemma}
\newtheorem{prop}[thmtool]{Proposition}
\newtheorem{clm}[thmtool]{Claim}
\newtheorem{conjecture}{Conjecture}
\newtheorem{problem}{Problem}
\newcommand{\Proof}{\noindent{\bf Proof.}\ \ }
\newcommand{\Remarks}{\noindent{\bf Remarks:}\ \ }
\newcommand{\Remark}{\noindent{\bf Remark:}\ \ }
\newcommand{\whp}{{\bf whp}\ }
\newcommand{\prob}{probability}
\newcommand{\rn}{random}
\newcommand{\rv}{random variable}
\newcommand{\hpg}{hypergraph}
\newcommand{\hpgs}{hypergraphs}
\newcommand{\subhpg}{subhypergraph}
\newcommand{\subhpgs}{subhypergraphs}
\newcommand{\bH}{{\bf H}}
\newcommand{\cH}{{\cal H}}
\newcommand{\cT}{{\cal T}}
\newcommand{\cF}{{\cal F}}
\newcommand{\cD}{{\cal D}}
\newcommand{\cC}{{\cal C}}

\begin{abstract}
The classical result of Erd\H os and R\'enyi asserts that the random
graph $G(n,p)$ experiences sharp phase transition around
$p=\frac{1}{n}$ -- for any $\epsilon>0$ and
$p=\frac{1-\epsilon}{n}$, all connected components of $G(n,p)$ are
typically of size $O_{\epsilon}(\log n)$, while for
$p=\frac{1+\epsilon}{n}$,  with high probability there exists a
connected component of size linear in $n$. We provide a very simple
proof of this fundamental result; in fact, we prove that in the
supercritical regime $p=\frac{1+\epsilon}{n}$, the random graph
$G(n,p)$ contains typically a path of linear length. We also discuss
applications of our technique to other random graph models and to
positional games.
\end{abstract}

\section{Introduction}
In their groundbreaking paper \cite{ER60} from 1960, Paul Erd\H os
and Alfr\'ed R\'enyi made the following fundamental discovery: the
random graph $G(n,p)$ undergoes a remarkable phase transition around
the edge probability $p(n)=\frac{1}{n}$. For any constant
$\epsilon>0$, if $p=\frac{1-\epsilon}{n}$, then $G(n,p)$ has {\bf
whp}\footnote{We say that an event  ${\cal E}_n$ occurs with high
probability, or \whp for brevity, in the probability space $G(n,p)$
if $\lim_{n\rightarrow\infty}Pr[G\sim G(n,p)\in {\cal E}_n]=1$.} all
connected components of size at most logarithmic in $n$, while for
$p=\frac{1+\epsilon}{n}$ \whp a connected component of linear size,
usually called the giant component, emerges in $G(n,p)$ (they also
showed that \whp there is a unique linear sized component). The
Erd\H os-R\'enyi paper, which launched the modern theory of random
graphs, has had enormous influence on the development of the field
and is generally considered to be a single most important paper in
Probabilistic Combinatorics, if not in all of Combinatorics.

There are now several proofs available for this result. Erd\H os and
R\'enyi (who actually worked in the model $G(n,m)$ of random graphs)
used counting arguments. Some of later proofs relied on the
machinery of branching processes. As one can expect for a result of
this magnitude of importance, there have been countless
ramifications and extensions proven over the years, and by now the
evolution of random graphs is very well understood. We refer the
reader to the standard sources in the theory of random graphs
\cite{JLR}, \cite{Bol-book} for a detailed account.

In 1981, Ajtai, Koml\'os and Szemer\'edi proved \cite{AKS81} that in
the supercritical regime $p=\frac{1+\epsilon}{n}$, not only the
random graph $G(n,p)$ contains \whp a linear sized connected
component, but it typically has a path of length linear in $n$.

The purpose of this note is to present a very simple and
self-contained proof of the Erd\H os-R\'enyi result, as well of the
result of Ajtai, Koml\'os and Szemer\'edi. We do not strive to
derive the best possible absolute constants, aiming rather for
simplicity.

Our notation is fairly standard. We set $N=\binom{n}{2}$. Floor and
ceiling signs will be systematically omitted for the sake of clarity
of presentation.

\section{Main result}

Our argument will utilize the notion of the Depth First Search
(DFS). This is a well known graph exploration algorithm, and we thus
will describe it rather briefly.

Recall that the DFS (Depth First Search) is a graph search algorithm
that visits all vertices of a (directed or undirected) graph
$G=(V,E)$ as follows. It maintains three sets of vertices, letting
$S$ be the set of vertices whose exploration is complete, $T$ be the
set of unvisited vertices, and $U=V\setminus(S\cup T)$, where the
vertices of $U$ are kept in a stack (the last in, first out data
structure). It is also assumed that some order $\sigma$ on the
vertices of $G$ is fixed, and the algorithm prioritizes vertices
according to $\sigma$. The algorithm starts with $S=U=\emptyset$ and
$T=V$, and runs till $U\cup T=\emptyset$. At each round of the
algorithm, if the set $U$ is non-empty, the algorithm queries $T$
for neighbors of the last vertex $v$ that has been added to $U$,
scanning $T$ according to $\sigma$. If $v$ has a neighbor $u$  in
$T$, the algorithm deletes $u$ from $T$ and inserts it into $U$. If
$v$ does not have a neighbor in $T$, then $v$ is popped out of $U$
and is moved to $S$. If $U$ is empty, the algorithm chooses the
first vertex of $T$ according to $\sigma$, deletes it from $T$ and
pushes it into $U$. In order to complete the exploration of the
graph, whenever the sets $U$ and $T$ have both become empty (at this
stage the connected component structure of $G$ has already been
revealed), we make the algorithm query all remaining pairs of
vertices in $S=V$, not queried before.

Observe that the DFS algorithm starts revealing a connected
component $C$ of $G$ at the moment the first vertex of $C$ gets into
(empty beforehand) $U$ and completes discovering all of $C$ when $U$
becomes empty again. We call a period of time between two
consecutive emptyings of $U$ an {\em epoch}, each epoch corresponds
to one connected component of $G$.

The following properties of the DFS algorithm will be relevant to
us:
\begin{itemize}
\item at each round of the algorithm one vertex moves, either from
$T$ to $U$, or from $U$ to $S$;
\item at any stage of the algorithm, it has been revealed already
that the graph $G$ has no edges between the current set $S$ and the
current set $T$;
\item the set $U$ always spans a path (indeed, when a vertex $u$ is
added to $U$, it happens because $u$ is a neighbor of the last
vertex $v$ in $U$; thus, $u$ augments the path spanned by $U$, of
which $v$ is the last vertex).
\end{itemize}

We will run the DFS on a random input $G\sim G(n,p)$, fixing the
order $\sigma$ on $V(G)=[n]$ to be the identity permutation. When
the DFS algorithm is fed with a sequence of i.i.d. Bernoulli($p$)
random variables $\bar{X}=(X_i)_{i=1}^N$, so that is gets its $i$-th
query answered positively if $X_i=1$ and answered negatively
otherwise, the so obtained graph is clearly distributed according to
$G(n,p)$. Thus, studying the component structure of $G$ can be
reduced to studying the properties of the random sequence $\bar{X}$.
In particular, observe crucially that as long as $T\ne\emptyset$,
every positive answer to a query results in a vertex being moved
from $T$ to $U$, and thus after $t$ queries and assuming
$T\ne\emptyset$ still, we have $|S\cup U|\ge \sum_{i=1}^t X_i$. (The
last inequality is strict in fact as the first vertex of each
connected component is moved from $T$ to $U$ ``for free", i.e.,
without need to get a positive answer to a query.) On the other
hand, since the addition of every vertex, but the first one in a
connected component, to $U$ is caused by a positive answer to a
query, we have at time $t$: $|U|\le 1+\sum_{i=1}^t X_i$.

The probabilistic part of our argument is provided by the following
quite simple lemma.

\begin{lemma}\label{le1}
Let $\epsilon>0$ be a small enough constant. Consider the sequence
$\bar{X}=(X_i)_{i=1}^N$ of i.i.d. Bernoulli random variables with
parameter $p$.
\begin{enumerate}
\item Let $p=\frac{1-\epsilon}{n}$.  Let $k=\frac{7}{\epsilon^2}\ln
n$. Then \whp there is no interval of length $kn$ in $[N]$, in which
at least $k$ of the random variables $X_i$ take value 1.
\item Let $p=\frac{1+\epsilon}{n}$. Let $N_0=\frac{\epsilon
n^2}{2}$. Then \whp $\left|\sum_{i=1}^{N_0} X_i
-\frac{\epsilon(1+\epsilon)n}{2}\right|\le n^{2/3}$.
\end{enumerate}
\end{lemma}

\Proof {1)} For a given interval $I$ of length $kn$ in $[N]$, the
sum $\sum_{i\in I} X_i$ is distributed binomially with parameters
$kn$ and $p$. Applying the standard Chernoff-type bound (see, e.g.,
Theorem A.1.11 of \cite{AS}) to the upper tail of $B(kn,p)$, and
then the union bound, we see that the probability of the existence
of an interval violating the assertion of the lemma is at most
$$
(N-k+1)Pr[B(kn,p)\ge k]< n^2\cdot
e^{-\frac{\epsilon^2}{3}(1-\epsilon)k} <n^2\cdot
e^{-\frac{\epsilon^2(1-\epsilon)}{3}\,\frac{7}{\epsilon^2}\ln n}=
o(1)\,,
$$
for small enough $\epsilon>0$.

\noindent{2)} The sum $\sum_{i=1}^{N_0}X_i$ is distributed
binomially with parameters $N_0$ and $p$. Hence, its expectation is
$N_0p=\frac{\epsilon n^2p}{2}=\frac{\epsilon(1+\epsilon)n}{2}$, and
its standard deviation is of order $n$. Applying the Chebyshev
inequality, we get the required estimate.

\medskip

Now we are ready to formulate and to prove our main result.

\begin{thm}\label{th1}
Let $\epsilon>0$ be a small enough constant. Let $G\sim G(n,p)$.
\begin{enumerate}
\item Let $p=\frac{1-\epsilon}{n}$. Then \whp all connected
components of $G$ are of size at most $\frac{7}{\epsilon^2}\ln n$.
\item Let $p=\frac{1+\epsilon}{n}$. Then \whp $G$ contains a path of
length at least $\frac{\epsilon^2 n}{5}$.
\end{enumerate}
\end{thm}
In both cases, we run the DFS algorithm on $G\sim G(n,p)$, and
assume that the sequence $\bar{X}=(X_i)_{i=1}^N$ of random
variables, defining the random graph $G\sim G(n,p)$ and guiding the
DFS algorithm, satisfies the corresponding part of Lemma \ref{le1}.

\medskip

\Proof 1) Assume to the contrary that $G$ contains a connected
component $C$ with more than $k=\frac{7}{\epsilon^2}\ln n$ vertices.
Let us look at the epoch of the DFS when $C$ was created. Consider
the moment inside this epoch when the algorithm has found the
$(k+1)$-st vertex of $C$ and is about to move it to $U$. Denote
$\Delta S=S\cap C$ at that moment. Then $|\Delta S\cup U|=k$, and
thus the algorithm got exactly $k$ positive  answers to its queries
to random variables $X_i$ during the epoch, with each positive
answer being responsible for revealing a new vertex of $C$, after
the first vertex of $C$ was put into $U$ in the beginning of the
epoch. At that moment during the epoch only pairs of edges touching
$\Delta S\cup U$ have been queried, and the number of such pairs is
therefore at most $\binom{k}{2}+k(n-k)<kn$. It thus follows that the
sequence $\bar{X}$ contains an interval of length at most $kn$ with
at least $k$ 1's inside -- a contradiction to Property 1 of Lemma
\ref{le1}.

2) Assume that the sequence $\bar{X}$ satisfies Property 2 of Lemma
\ref{le1}. We claim that after the first $N_0=\frac{\epsilon
n^2}{2}$ queries of the DFS algorithm, the set $U$ contains at least
$\frac{\epsilon^2 n}{5}$ vertices (with the contents of $U$ forming
a path of desired length at that moment). Observe first that
$|S|<\frac{n}{3}$ at time $N_0$. Indeed, if $|S|\ge\frac{n}{3}$,
then let us look at a moment $t$ where $|S|=\frac{n}{3}$ (such a
moment surely exists as vertices flow to $S$ one by one). At that
moment $|U|\le 1+\sum_{i=1}^{t}X_i<\frac{n}{3}$ by Property 2 of
Lemma \ref{le1}. Then $|T|=n-|S|-|U|\ge \frac{n}{3}$, and the
algorithm has examined all $|S|\cdot|T|\ge \frac{n^2}{9}>N_0$ pairs
between $S$ and $T$ (and found them to be non-edges) -- a
contradiction. Let us return to time $N_0$. If $|S|<\frac{n}{3}$ and
$|U|<\frac{\epsilon^2 n}{5}$ then, we have $T\ne\emptyset$. This
means in particular that the algorithm is still revealing the
connected components of $G$, and each positive answer it got
resulted in moving a vertex from $T$ to $U$ (some of these vertices
may have already moved further from $U$ to $S$). By Property 2 of
Lemma \ref{le1} the number of positive answers at that point is at
least $\frac{\epsilon(1+\epsilon)n}{2}-n^{2/3}$. Hence we have
$|S\cup U|\ge \frac{\epsilon(1+\epsilon)n}{2}-n^{2/3}$. If
$|U|\le\frac{\epsilon^2n}{5}$, then $|S|\ge \frac{\epsilon
n}{2}+\frac{3\epsilon^2n}{10}-n^{2/3}$. All $|S||T|\ge
|S|\left(n-|S|-\frac{\epsilon^2n}{5}\right)$ pairs between $S$ and
$T$ have been probed by the algorithm (and answered in the
negative). We thus get:
\begin{eqnarray*}
\frac{\epsilon n^2}{2}&=&N_0\ge
|S|\left(n-|S|-\frac{\epsilon^2n}{5}\right)\ge \left(\frac{\epsilon
n}{2}+\frac{3\epsilon^2n}{10}-n^{2/3}\right) \left(n-\frac{\epsilon
n}{2}-\frac{\epsilon^2n}{2}+n^{2/3}\right)\\
&=&\frac{\epsilon n^2}{2}+\frac{\epsilon^2n^2}{20}-O(\epsilon^3)n^2
> \frac{\epsilon n^2}{2}
\end{eqnarray*}
(we used the assumption $|S|<\frac{n}{3})$, and this is obviously a
contradiction, completing the proof.

\section{Discussion}

\noindent{\bf 1.} Observe that using a Chernoff-type bound for the
tales of the binomial random variable instead of the Chebyshev
inequality would allow to claim in the second part of Lemma
\ref{le1} that the sum $\sum_{i=1}^{N_0} X_i$ is close to
$\frac{\epsilon(1+\epsilon)n}{2}$ with probability exponentially
close to 1. This would show in turn, employing the argument of
Theorem \ref{th1}, that $G(n,p)$ with $p=\frac{1+\epsilon}{n}$
contains a path of length linear in $n$ with exponentially high
probability, namely, with probability $1-\exp\{-c(\epsilon)n\}$.

\medskip

 \noindent{\bf 2.} The dependencies on $\epsilon$ in
both parts of Theorem \ref{th1} are of the correct order of
magnitude -- for $p=\frac{1-\epsilon}{n}$ a largest connected
component of $G(n,p)$ is known to be \whp of size
$\Theta(\epsilon^{-2})\log n$ (see, e.g., Cors. 5.8 and 5.11 of
\cite{Bol-book}), while for $p=\frac{1+\epsilon}{n}$ a longest cycle
of $G(n,p)$ is \whp of length $\Theta(\epsilon^2)n$ (see, e.g., Th.
5.17 of \cite{JLR}); the standard trick of sprinkling further random
edges with edge probability $p'=o(n^{-1})$ shows that if $G(n,p)$
contains \whp a path of length $\alpha n$ for some constant
$\alpha>0$, then $G(n,p+p')$ contains \whp a cycle of length at
least $(\alpha-o(1))n$. Note also that although we stated our result
in Theorem \ref{th1} for a constant $\epsilon>0$, our argument is in
fact valid for $\epsilon=\epsilon(n)\rightarrow 0$ as well, with a
bit more careful treatment of the error terms in our proofs.
Actually, we can take $\epsilon(n)$ to be as low as
$\epsilon \gg n^{-1/3}\log^{1/3} n$ in our arguments (including
the theorem in the next remark) -- which nearly borders the critical
window $\epsilon=\Theta(n^{-1/3})$.

\medskip

\noindent{\bf 3.} The giant component itself in the regime
$p=\frac{1+\epsilon}{n}$, $\epsilon>0$ a constant, is known to be
substantially larger typically than a longest path -- it has \whp
$\Theta(\epsilon)n$ vertices (see, e.g., Th. 5.4 of \cite{JLR}).
Using very similar techniques, we can show the probable existence of
a connected component of size $\Omega(\epsilon)n$ in this range, as
given by the following theorem.

\begin{thm}\label{thm2n}
Let $p=\frac{1+\epsilon}{n}$, for $\epsilon>0$ a small enough
constant. Let $G\sim G(n,p)$. Then \whp $G$ has a connected
component with at least $\frac{\epsilon n}{2}$ vertices.
\end{thm}

\Proof The proof is quite similar to that of Theorem \ref{th1}, and
therefore we will allow ourselves to be rather concise. Here too we
run the DFS algorithm on $G\sim G(n,p)$ and feed it with a sequence
$\bar{X}$ of i.i.d. Bernoulli($p$) random variables
$\bar{X}=(X_i)_{i=1}^N$. Denote as before $N_0=\frac{\epsilon
n^2}{2}$. We will need the following typical properties of the
sequence $\bar{X}$, slightly generalizing those stated in Part 2. of
Lemma \ref{le1} and provable using the same Chernoff-type estimates:
\begin{enumerate}
\item $\sum_{i=1}^{n^{7/4}}X_i\le n^{5/6}$;
\item For every $n^{7/4}\le t\le N_0$,
$\left|\sum_{i=1}^t X_i-(1+\epsilon)\frac{t}{n}\right|\le n^{2/3}$.
\end{enumerate}

Let us assume now that the sequence $\bar{X}$ satisfies the above
stated properties. We claim that after the first $N_0$ queries of
the DFS algorithm, we are in the midst of revealing a connected
component whose size is at least $\frac{\epsilon n}{2}$. Just as in
the proof of Theorem \ref{th1} we have that $|S|<\frac{n}{3}$ at
time $N_0$, and $T$ is still non-empty. It follows that at any
moment $n^{7/4}\le t\le N_0$ we have: $|S\cup U|\ge
(1+\epsilon)\frac{t}{n}-n^{2/3}$. If at some moment $t$ in this
interval the set $U$ becomes empty, the algorithm has asked all
queries between the set $S$ and its complement $T=[n]-S$, implying:
\begin{eqnarray*}
t&\ge& |S|(n-|S|)\ge \left((1+\epsilon)\frac{t}{n}-n^{2/3}\right)
\left(n-(1+\epsilon)\frac{t}{n}+n^{2/3}\right) \ge
(1+\epsilon)t-(1+\epsilon)^2\frac{t^2}{n^2}-2n^{5/3}\\
&\ge& (1+\epsilon)\left(1-(1+\epsilon)\frac{\epsilon}{2}\right)t
-2n^{5/3}=
(1+\epsilon)\left(1-\frac{\epsilon}{2}-\frac{\epsilon^2}{2}\right)t-2n^{5/3}>t
\end{eqnarray*}
-- a contradiction, for small enough $\epsilon>0$. (We used
$|S|<\frac{n}{3}$ in the above estimate.) Hence $U$ is never empty
in the interval $[n^{7/4},N_0]$. It follows that all vertices added
to $U$ during this interval (of which some may have migrated further
to $S$) are in the same connected component, and their number is, by
the properties of $\bar{X}$ stated above,
$$
\sum_{i=n^{7/4}}^{N_0}X_i\ge
(1+\epsilon)\frac{N_0}{n}-n^{2/3}-n^{5/6} \ge
(1+\epsilon)\frac{\epsilon n}{2}-2n^{5/6}\ge \frac{\epsilon n}{2}\,.
$$
All these vertices belong to the same connected component -- whose
size is then at least $\frac{\epsilon n}{2}$, completing the proof.

\medskip

\noindent{\bf 4.} As we have already mentioned, the DFS algorithm is
applicable equally well to directed graphs. Hence essentially the
same argument as above, with obvious minor changes, can be applied
to the model $D(n,p)$ of random digraphs. In this model, the vertex
set is $[n]$, and each of the $n(n-1)$ ordered pairs $(i,j)$, $1\le
i\ne j\le n$, is a directed edge of $D\sim D(n,p)$ with probability
$p=p(n)$ and independently from other pairs. In particular we can
obtain the following theorem:
\begin{thm}\label{thm2}
Let $p=\frac{1+\epsilon}{n}$, for $\epsilon>0$ constant. Then the
random digraph $D(n,p)$ has \whp a directed path and a directed
cycle of length $\Theta(\epsilon^2)n$.
\end{thm}
This recovers the classical result of Karp \cite{Kar90} for the
model $D(n,p)$.

\medskip

\noindent{\bf 5.} The technique of Theorem \ref{th1} can be applied
to further models of random graphs and digraphs. One immediate
application is to random subgraphs of graphs of large minimum
degree. We have the following theorem.

\begin{thm}\label{thm-mindeg}
Let $G$ be a finite graph with minimum degree at least $n$. Let
$p=\frac{1+\epsilon}{n}$, for $\epsilon>0$ constant. Form a random
subgraph $G_p$ of $G$ by including every edge of $G$ into $G_p$
independently and with probability $p$. Then \whp $G_p$ has a path
of length at least $\frac{\epsilon^2n}{5}$.
\end{thm}

The proof is essentially identical to that of Theorem \ref{th1}. We
run the DFS process on $G_p$ and feed it with a sequence $\bar{X}$
of i.i.d. Bernoulli($p$) random variables $\bar{X}=(X_i)_{i=1}^N$,
where $N=|E(G)|$. For the proof, we need only to notice that at any
time the number of edges of $G$ between $S$ and $T$ can be estimated
from below by $|S|(\delta(G)-|S|-|U|)\ge|S|(n-|S|-|U|)$, the rest of
the proof is the same. Notice that getting a long cycle appears to
be a much more challenging task in this setting -- the base graph
$G$ can be of girth (much) larger than $n$, and therefore sprinkling
does not necessarily help (immediately) to turn a long path into a
long cycle \whp.

\medskip

\noindent{\bf 6.} Another example of applying our technique is
random subgraphs of pseudo-random graphs. Let $G$ be an
$(n,d,\lambda)$-graph (a $d$-regular graph on $n$ vertices, in which
all eigenvalues of the adjacency matrix, but the first one, are at
most $\lambda$ in their absolute values -- see, e.g. \cite{KS06} for
a thorough discussion of this notion). It is well known that
requiring $\lambda\ll d$ is enough to guarantee many pseudo-random
properties of such a graph. The model of taking a random subgraph
$G_p$ of an $(n,d,\lambda)$-graph $G$ has been considered by Frieze,
Krivelevich and Martin in \cite{FKM04}. It is proven in \cite{FKM04}
that, assuming $\lambda\ll d$, for $p=\frac{1+\epsilon}{d}$ the
random subgraph $G_p$ of an an $(n,d,\lambda)$-graph $G$ has \whp
the unique connected component of size linear in $n$. We can apply
the technique of Theorem \ref{th1} to prove the following:
\begin{thm}\label{thm3}
Let $G$ be an $(n,d,\lambda)$-graph with $\lambda=o(d)$. Let
$p=\frac{1+\epsilon}{d}$, for $\epsilon>0$ constant. Then the random
subgraph $G_p$ contains \whp\ a path of length
$\Theta(\epsilon^2)n$.
\end{thm}
Here is a very brief sketch of the proof. We run the DFS algorithm
on $G_p$ till it queries $\frac{\epsilon dn}{2}$ edges of $G$.

Similarly to Lemma \ref{le1}, it gets \whp about
$\frac{\epsilon(1+\epsilon)n}{2}$ positive answers during this
period, when fed with a string of i.i.d.
Bernoulli$\left(\frac{1+\epsilon}{d}\right)$ random variables. In
order for the proof analogous to that of Theorem \ref{th1} to go
through, one only needs to be able to control the number of edges
between any two linear sized vertex subsets $S,T$ in $G$. Such a
control is indeed available for $(n,d,\lambda)$-graphs -- it is
known that if $G$ is an $(n,d,\lambda)$-graph, then for any two
vertex subsets $S,T\subseteq V(G)$ the number $e_G(S,T)$ of edges of
$G$ with one endpoint in $S$ and another in $T$ satisfies:
$$
\left|e_G(S,T)-\frac{d}{n}|S|\,|T|\right|\le \lambda\sqrt{|S||T|}
$$
(see, e.g. Corollary 9.2.5 of \cite{AS} or Theorem 2.11 of
\cite{KS06}). Assuming $\lambda\ll d$ is enough therefore to
guarantee that $e_G(S,T)=(1+o(1))\frac{d}{n}|S|||T|$ in such a
graph, and the proof for the random subgraph proceeds as in Theorem
\ref{th1}. Here too sprinkling helps to turn a long path into a long
cycle \whp -- we first get \whp a linearly long path and then argue
that due to the above estimate on the edge distribution of $G$ there
are $\Theta(dn)$ edges between the prefix and the suffix of the
path, and one of them will \whp fall into a sprinkled graph, thus
closing a long cycle.

\medskip

\noindent{\bf 7.} Yet another application of our proof strategy is
to positional games. The following game ${\cal L}(n,b)$ was
considered by Bednarska and \L uczak in \cite{BL01}. The game is
played between two players, Maker and Breaker, alternately claiming
1 and $b$ edges, respectively, of the complete graph $K_n$ on $n$
vertices, till all edges of $K_n$ have been claimed by either of the
players. Maker's goal is to maximize the number of vertices in a
largest connected component in her graph by the end of the game,
Breakers aims to make it as small as possible. Bednarska and \L
uczak discovered the following phase transition phenomenon,
obviously reminiscent of the Erd\H os-R\'enyi phase transition in
random graphs. Let $\epsilon>0$ be a constant. If $b=(1+\epsilon)n$
then Breaker has a strategy to keep all of Maker's connected
components of size $O(1/\epsilon)$. On the other hand, if
$b=(1-\epsilon)n$, then Maker has a strategy to create a connected
component of size $\Theta(\epsilon)n$. We can prove the following
result.
\begin{thm}\label{thm4}
Let $\epsilon>0$. Then in the game ${\cal L}(n,b)$ with
$b=(1-\epsilon)n$, Maker has a strategy to create a path of length
$\Theta(\epsilon^2)n$.
\end{thm}
The winning strategy of Maker and the proof of its validity are
fairly similar to the proof of Theorem \ref{th1}. Maker maintains
three sets $S,U,T$ partitioning $[n]$, starting with $S=\emptyset$,
and $U$ being an arbitrary vertex from $[n]$. She makes sure that
the set $U$ always spans a path of her edges at any stage of the
game. At each Maker's turn, she finds the last vertex $v$ along the
path in $U$ for which there exists an unclaimed edge $(v,u)$ with
$u\in T$, shifts all further vertices after $v$ along $U$ into $S$
and claims the edge $(v,u)$, moving $u$ from $T$ to $U$. If no such
vertex is available along the current path in $U$, Maker moves all
of its vertices into $S$, loads $U$ with an arbitrary vertex $u$
from $T$ and then proceeds as described before. One can observe
that, similarly to the analysis of the DFS algorithm, at any stage
of the game all edges between the current set $S$ and the current
set $T$ have been claimed by Breaker. Now, look at the situation in
the game after $\frac{\epsilon n}{2}$ rounds. At that point $|S\cup
U|\ge \frac{\epsilon n}{2}$. If one has $|U|\le
\frac{\epsilon^2n}{5}$, then all
$$
|S|\,|T|\ge \left(\frac{\epsilon
n}{2}-\frac{\epsilon^2n}{5}\right)\,\left(n-\frac{\epsilon
n}{2}\right)> \frac{\epsilon n}{2}(1-\epsilon)n
$$
edges between $S$ and $T$ have been claimed by Breaker -- a
contradiction, for small enough $\epsilon>0$. The situation with
making a cycle is quite different here -- it has been shown by
Bednarska and Pikhurko \cite{BP05} that if $b=b(n)$ is such that
Maker completes the game with at most $n-1$ edges, then Breaker has
a strategy to force Maker to end up with a tree; thus $b\ge
(1+o(1))n/2$ is required for Maker to create a cycle of any length.

\medskip

\noindent{\bf 8.} Some of the idea utilized in this paper have
already been applied before. In particular, the DFS algorithm has
been used by Ben-Eliezer and the authors in \cite{BKS12} to prove
the following statement: if in a graph $G$ on $n$ vertices there is
an edge between every pair of disjoint vertex subsets of size $k$,
then $G$ contains a path of length $n-2k+1$. This deterministic
statement implies readily that $G(n,p)$ with $p=c/n$ contains \whp a
path of length $(1-\alpha(c))n$, where $\alpha(c)\rightarrow 0$ as
$c\rightarrow\infty$. Also, Benjamini and Schramm \cite{BS96} used
the idea of coupling a graph search algorithm with a sequence
$\bar{X}$ of random bits, serving as answers to the algorithm's
queries, to derive some results about percolation in expanding
graphs.


\end{document}